\documentclass[11pt]{amsart}
\usepackage{amsmath,amssymb}
\newtheorem*{theorem}{Theorem}
\newtheorem{proposition}{Proposition}
\newtheorem{corollary}[proposition]{Corollary}

\begin{document}

\title{Blow-up phenomena for the Yamabe equation II}
\author{Simon Brendle and Fernando C. Marques}

\maketitle

\section{Introduction}

Let $(M,g)$ be a compact Riemannian manifold of dimension $n \geq 3$.
The Yamabe problem is concerned with finding metrics of constant scalar
curvature in the conformal class of $g$. This problem leads to a semi-linear elliptic PDE for the conformal factor. More precisely, a conformal metric of the form $u^{\frac{4}{n-2}} \, g$ has constant scalar curvature $c$ if and only if
\begin{equation}
\label{yamabe.pde}
\frac{4(n-1)}{n-2} \, \Delta_g u - R_g \, u + c \, u^{\frac{n+2}{n-2}} = 0,
\end{equation}
where $\Delta_g$ is the Laplace operator with respect to $g$ and $R_g$ denotes the scalar curvature of $g$. Every solution of (\ref{yamabe.pde}) is a critical point of the functional
\begin{equation}
\label{yamabe.functional}
E_g(u) = \frac{\int_M \big ( \frac{4(n-1)}{n-2} \, |du|_g^2 + R_g \, u^2
\big ) \, dvol_g}{\big ( \int_M
u^{\frac{2n}{n-2}} \, dvol_g \big )^{\frac{n-2}{n}}}.
\end{equation}
In this paper, we address the question whether the set of all solutions to the Yamabe PDE is compact in the $C^2$-topology. It has been conjectured that this should be true unless $(M,g)$ is conformally equivalent to the round sphere (see \cite{Schoen1},\cite{Schoen2},\cite{Schoen3}). The case of the round sphere $S^n$ is special in that (\ref{yamabe.pde}) is invariant under the action of the conformal group on $S^n$, which is non-compact. It follows from a theorem of Obata \cite{Obata} that every solution of the Yamabe PDE on $S^n$ is minimizing, and the space of all solutions to the Yamabe PDE on $S^n$
can be identified with the unit ball $B^{n+1}$. Note that the round sphere is the only compact manifold for which the set of minimizing solutions is non-compact.

The Compactness Conjecture has been verified in low dimensions and in the locally conformally flat case. If $(M,g)$ is locally conformally flat, compactness follows from work of R.~Schoen \cite{Schoen1},\cite{Schoen2}. Moreover, Schoen proposed a strategy, based on the Pohozaev identity, for 
proving the conjecture in the non-locally conformally flat case. In \cite{Li-Zhu}, Y.Y.~Li and M.~Zhu followed this strategy to prove compactness in dimension 
$3$. O.~Druet \cite{Druet} proved the conjecture in dimensions $4$ and $5$.

The case $n \geq 6$ is more subtle, and requires a careful analysis of
the local properties of the background metric $g$ near a blow-up point. 
The Compactness Conjecture is closely related to the Weyl Vanishing 
Conjecture, which  asserts that the Weyl tensor should vanish to an order 
greater than $\frac{n-6}{2}$ at a blow-up point (see \cite{Schoen3}). The 
Weyl Vanishing Conjecture has been verified in dimensions $6$ and $7$ by
F.~Marques \cite{Marques} and, independently, by Y.Y.~Li and L.~Zhang
\cite{Li-Zhang1}. Using these results and the positive mass theorem, these 
authors were able to prove compactness for $n \leq 7$. Moreover, Li and Zhang 
showed that compactness holds in all dimensions provided that $|W_g(p)| 
+ |\nabla W_g(p)| > 0$ for all $p \in M$. In dimensions $10$ and $11$, it 
is sufficient to assume that $|W_g(p)| + |\nabla W_g(p)| + |\nabla^2 W_g(p)| > 0$ for all $p \in M$
(see \cite{Li-Zhang2}).

Very recently, M.~Khuri, F.~Marques and R.~Schoen \cite{Khuri-Marques-Schoen} 
proved the Weyl Vanishing Conjecture up to dimension $24$. This result, 
combined with the positive mass theorem, implies the Compactness Conjecture for
those dimensions. After proving sharp pointwise estimates, they reduce
these questions to showing a certain quadratic
form is positive definite. It turns out the quadratic form has negative
eigenvalues if $n \geq 25$.

In a recent paper \cite{Brendle1}, it was shown that the Compactness Conjecture 
fails for $n \geq 52$. More precisely, given any integer $n \geq 52$, there 
exists a smooth Riemannian metric $g$ on $S^n$ such that set of constant scalar 
curvature metrics in the conformal class of $g$ is non-compact. Moreover, the 
blowing-up sequences obtained in \cite{Brendle1} form exactly one bubble. The 
construction relies on a gluing procedure based on some local model metric. 
These local models are directions in which the quadratic form of 
\cite{Khuri-Marques-Schoen} is negative definite. We refer to \cite{Brendle2} for a 
survey of this and related results.

In the present paper, we extend these counterexamples to the dimensions
$25 \leq n \leq 51$. Our main theorem is:

\begin{theorem}
\label{main.theorem}
Assume that $25 \leq n \leq 51$. Then there exists a Riemannian metric
$g$ on $S^n$ (of class $C^\infty$) and a sequence of positive functions
$v_\nu \in C^\infty(S^n)$ ($\nu \in \mathbb{N}$) with the following
properties:
\begin{itemize}
\item[(i)] $g$ is not conformally flat
\item[(ii)] $v_\nu$ is a solution of the Yamabe PDE (\ref{yamabe.pde})
for all $\nu \in \mathbb{N}$
\item[(iii)] $E_g(v_\nu) < Y(S^n)$ for all $\nu \in \mathbb{N}$, and
$E_g(v_\nu) \to Y(S^n)$ as $\nu \to \infty$
\item[(iv)] $\sup_{S^n} v_\nu \to \infty$ as $\nu \to \infty$
\end{itemize}
(Here, $Y(S^n)$ denotes the Yamabe energy of the round metric on $S^n$.)
\end{theorem}

We note that O.~Druet and E.~Hebey \cite{Druet-Hebey1} have constructed 
blow-up examples for perturbations of (\ref{yamabe.pde}) (see also \cite{Druet-Hebey2}).

In Section 2, we describe how the problem can be reduced to finding critical points of a certain function $\mathcal{F}_g(\xi,\varepsilon)$, where $\xi$ is a vector in $\mathbb{R}^n$ and $\varepsilon$ is a positive real number. This idea has been used by many authors (see, e.g., \cite{Ambrosetti}, \cite{Ambrosetti-Malchiodi}, \cite{Berti-Malchiodi}, \cite{Brendle1}). In Section 3, we show that the function
$\mathcal{F}_g(\xi,\varepsilon)$ can be approximated by an auxiliary
function $F(\xi,\varepsilon)$. In Section 4, we prove that the function
$F(\xi,\varepsilon)$ has a critical point, which is a strict local
minimum. Finally, in Section 5, we use a perturbation argument to
construct critical points of the function
$\mathcal{F}_g(\xi,\varepsilon)$. From this the non-compactness result
follows.  

The authors would like to thank Professor Richard Schoen for constant support and encouragement. The first author was supported by the Alfred P. Sloan foundation and by the National Science Foundation under grant
DMS-0605223. The second author was supported by CNPq-Brazil, FAPERJ and the Stanford Mathematics Department. \\

\section{Lyapunov-Schmidt reduction}

In this section, we collect some basic results established in
\cite{Brendle1}. Let
\[\mathcal{E} = \bigg \{ w \in L^{\frac{2n}{n-2}}(\mathbb{R}^n) \cap
W_{loc}^{1,2}(\mathbb{R}^n): \int_{\mathbb{R}^n} |dw|^2 < \infty \bigg
\}.\] By Sobolev's inequality, there exists a constant $K$, depending
only on $n$, such that
\[\bigg ( \int_{\mathbb{R}^n} |w|^{\frac{2n}{n-2}} \bigg
)^{\frac{n-2}{n}} \leq K \, \int_{\mathbb{R}^n} |dw|^2\]
for all $w \in \mathcal{E}$. We define a norm on $\mathcal{E}$ by
$\|w\|_{\mathcal{E}}^2 = \int_{\mathbb{R}^n} |dw|^2$. It is easy to see
that $\mathcal{E}$, equipped with this norm, is complete. \\

Given any pair $(\xi,\varepsilon) \in \mathbb{R}^n \times (0,\infty)$,
we define a function $u_{(\xi,\varepsilon)}:
\mathbb{R}^n \to \mathbb{R}$ by
\[u_{(\xi,\varepsilon)}(x) = \Big ( \frac{\varepsilon}{\varepsilon^2 +
|x - \xi|^2} \Big )^{\frac{n-2}{2}}.\]
The function $u_{(\xi,\varepsilon)}$ satisfies the elliptic PDE
\[\Delta u_{(\xi,\varepsilon)} + n(n-2) \,
u_{(\xi,\varepsilon)}^{\frac{n+2}{n-2}} = 0.\]
It is well known that
\[\int_{\mathbb{R}^n} u_{(\xi,\varepsilon)}^{\frac{2n}{n-2}} = \Big (
\frac{Y(S^n)}{4n(n-1)} \Big )^{\frac{n}{2}}\]
for all $(\xi,\varepsilon) \in \mathbb{R}^n \times (0,\infty)$. We next
define
\[\varphi_{(\xi,\varepsilon,0)}(x) =
\Big ( \frac{\varepsilon}{\varepsilon^2 + |x - \xi|^2} \Big
)^{\frac{n+2}{2}} \, \frac{\varepsilon^2 - |x - \xi|^2}{\varepsilon^2 +
|x - \xi|^2}\]
and
\[\varphi_{(\xi,\varepsilon,k)}(x) =
\Big ( \frac{\varepsilon}{\varepsilon^2 + |x - \xi|^2} \Big
)^{\frac{n+2}{2}} \, \frac{2\varepsilon \, (x_k - \xi_k)}{\varepsilon^2
+ |x - \xi|^2}\]
for $k = 1,\hdots,n$. Finally, we define a closed subspace
$\mathcal{E}_{(\xi,\varepsilon)} \subset \mathcal{E}$ by
\[\mathcal{E}_{(\xi,\varepsilon)} = \bigg \{ w \in \mathcal{E}:
\int_{\mathbb{R}^n} \varphi_{(\xi,\varepsilon,k)} \, w = 0 \quad
\text{for $k = 0,1,\hdots,n$} \bigg \}.\]
Clearly, $u_{(\xi,\varepsilon)} \in \mathcal{E}_{(\xi,\varepsilon)}$.

\begin{proposition}
\label{linearized.operator}
Consider a Riemannian metric on $\mathbb{R}^n$ of the form $g(x) =
\exp(h(x))$, where $h(x)$ is a trace-free symmetric two-tensor on
$\mathbb{R}^n$ satisfying  $h(x) = 0$ for $|x| \geq 1$. There exists a
positive constant $\alpha_0 \leq 1$, depending only on $n$, with the
following significance: if $|h(x)| + |\partial h(x)| + |\partial^2 h(x)|
\leq \alpha_0$ for all $x \in \mathbb{R}^n$, then, given any pair
$(\xi,\varepsilon) \in \mathbb{R}^n \times (0,\infty)$ and any function
$f \in L^{\frac{2n}{n+2}}(\mathbb{R}^n)$, there exists a unique function
$w=G_{(\xi,\varepsilon)}(f) \in \mathcal{E}_{(\xi,\varepsilon)}$ such that
\[\int_{\mathbb{R}^n} \Big ( \langle dw,d\psi \rangle_g
+ \frac{n-2}{4(n-1)} \, R_g \, w \, \psi - n(n+2) \,
u_{(\xi,\varepsilon)}^{\frac{4}{n-2}} \, w \, \psi \Big ) =
\int_{\mathbb{R}^n} f \, \psi\]
for all test functions $\psi \in \mathcal{E}_{(\xi,\varepsilon)}$.
Moreover, we have $\|w\|_{\mathcal{E}} \leq C \,
\|f\|_{L^{\frac{2n}{n+2}}(\mathbb{R}^n)}$, where $C$ is a constant that
depends only on $n$.
\end{proposition}

\begin{proposition}
\label{fixed.point.argument}
Consider a Riemannian metric on $\mathbb{R}^n$ of the form $g(x) =
\exp(h(x))$, where $h(x)$ is a trace-free symmetric two-tensor on
$\mathbb{R}^n$ satisfying $h(x) = 0$ for $|x| \geq 1$. Moreover, let
$(\xi,\varepsilon) \in \mathbb{R}^n \times (0,\infty)$. There exists a
positive constant $\alpha_1 \leq \alpha_0$, depending only on $n$, with
the following significance: if $|h(x)| + |\partial h(x)| + |\partial^2
h(x)| \leq \alpha_1$ for all $x \in \mathbb{R}^n$, then there exists a
function $v_{(\xi,\varepsilon)} \in \mathcal{E}$ such that
$v_{(\xi,\varepsilon)} - u_{(\xi,\varepsilon)} \in
\mathcal{E}_{(\xi,\varepsilon)}$
and
\[\int_{\mathbb{R}^n} \Big ( \langle dv_{(\xi,\varepsilon)},d\psi
\rangle_g + \frac{n-2}{4(n-1)} \, R_g \, v_{(\xi,\varepsilon)} \, \psi -
n(n-2) \, |v_{(\xi,\varepsilon)}|^{\frac{4}{n-2}} \,
v_{(\xi,\varepsilon)} \, \psi \Big ) = 0\] for all test functions $\psi
\in \mathcal{E}_{(\xi,\varepsilon)}$. Moreover, we have the estimate
\begin{align*}
&\|v_{(\xi,\varepsilon)} - u_{(\xi,\varepsilon)}\|_{\mathcal{E}} \\
&\leq C \, \Big \| \Delta_g u_{(\xi,\varepsilon)} - \frac{n-2}{4(n-1)}
\, R_g \, u_{(\xi,\varepsilon)} + n(n-2) \,
u_{(\xi,\varepsilon)}^{\frac{n+2}{n-2}} \Big
\|_{L^{\frac{2n}{n+2}}(\mathbb{R}^n)},
\end{align*}
where $C$ is a constant that depends only on $n$.
\end{proposition}

We next define a function $\mathcal{F}_g: \mathbb{R}^n \times (0,\infty)
\to \mathbb{R}$ by
\begin{align*}
\mathcal{F}_g(\xi,\varepsilon) &= \int_{\mathbb{R}^n}
\Big ( |dv_{(\xi,\varepsilon)}|_g^2 + \frac{n-2}{4(n-1)} \, R_g \,
v_{(\xi,\varepsilon)}^2 - (n-2)^2 \,
|v_{(\xi,\varepsilon)}|^{\frac{2n}{n-2}} \Big ) \\
&- 2(n-2) \, \Big ( \frac{Y(S^n)}{4n(n-1)} \Big )^{\frac{n}{2}}.
\end{align*}
If we choose $\alpha_1$ small enough, then we obtain the following
result: \\

\begin{proposition}
\label{reduction.to.a.finite.dimensional.problem}
The function $\mathcal{F}_g$ is continuously differentiable. Moreover,
if $(\bar{\xi},\bar{\varepsilon})$ is a critical point of the function
$\mathcal{F}_g$, then
the function $v_{(\bar{\xi},\bar{\varepsilon})}$ is a non-negative weak
solution of the equation
\[\Delta_g v_{(\bar{\xi},\bar{\varepsilon})} - \frac{n-2}{4(n-1)} \, R_g
\, v_{(\bar{\xi},\bar{\varepsilon})} + n(n-2) \,
v_{(\bar{\xi},\bar{\varepsilon})}^{\frac{n+2}{n-2}} = 0.\]
\end{proposition}


\section{An estimate for the energy of a ``bubble"}

Throughout this paper, we fix a real number $\tau$ and a multi-linear
form $W: \mathbb{R}^n \times \mathbb{R}^n \times \mathbb{R}^n \times
\mathbb{R}^n \to \mathbb{R}$. The number $\tau$ depends only on the
dimension $n$. The exact choice of $\tau$ will be postponed until
Section 4. We assume that $W_{ijkl}$ satisfies all the algebraic
properties of the Weyl tensor. Moreover, we assume that some components
of $W$ are non-zero, so that
\[\sum_{i,j,k,l=1}^n (W_{ijkl} + W_{ilkj})^2 > 0.\]
For abbreviation, we put
\[H_{ik}(x) = \sum_{p,q=1}^n W_{ipkq} \, x_p \, x_q\]
and
\[\overline{H}_{ik}(x) = f(|x|^2) \, H_{ik}(x),\]
where $f(s) = \tau + 5s - s^2 + \frac{1}{20} \, s^3$. It is easy to see
that $H_{ik}(x)$ is trace-free, $\sum_{i=1}^n x_i \, H_{ik}(x) = 0$, and
$\sum_{i=1}^n \partial_i H_{ik}(x) = 0$ for all $x \in \mathbb{R}^n$. \\

We consider a Riemannian metric of the form $g(x) = \exp(h(x))$, where
$h(x)$ is a trace-free symmetric two-tensor on $\mathbb{R}^n$ satisfying
$h(x) = 0$ for $|x| \geq 1$,
\[|h(x)| + |\partial h(x)| + |\partial^2 h(x)| \leq \alpha_1\]
for all $x \in \mathbb{R}^n$, and
\[h_{ik}(x) = \mu \, \lambda^6 \, f(\lambda^{-2} \, |x|^2) \, H_{ik}(x)\]
for $|x| \leq \rho$. We assume that the parameters $\lambda$, $\mu$, and
$\rho$ are chosen such that
$\mu \leq 1$ and $\lambda \leq \rho \leq 1$. Note that
$\sum_{i=1}^n x_i \, h_{ik}(x) = 0$ and $\sum_{i=1}^n \partial_i
h_{ik}(x) = 0$ for $|x| \leq \rho$.

Given any pair $(\xi,\varepsilon) \in \mathbb{R}^n \times (0,\infty)$,
there exists a unique function $v_{(\xi,\varepsilon)}$ such that $v_{(\xi,\varepsilon)}
- u_{(\xi,\varepsilon)} \in \mathcal{E}_{(\xi,\varepsilon)}$
and
\[\int_{\mathbb{R}^n} \Big ( \langle dv_{(\xi,\varepsilon)},d\psi
\rangle_g + \frac{n-2}{4(n-1)} \, R_g \, v_{(\xi,\varepsilon)} \, \psi -
n(n-2) \, |v_{(\xi,\varepsilon)}|^{\frac{4}{n-2}} \,
v_{(\xi,\varepsilon)} \, \psi \Big ) = 0\] for all test functions $\psi
\in \mathcal{E}_{(\xi,\varepsilon)}$ (see Proposition \ref{fixed.point.argument}).

For abbreviation, let
\[\Omega = \bigg \{ (\xi,\varepsilon) \in \mathbb{R}^n \times
\mathbb{R}: |\xi| < 1, \, \frac{1}{2} < \varepsilon < 2 \bigg \}.\]

The following result is proved in the Appendix A of \cite{Brendle1}. A similar formula is derived in \cite{Ambrosetti-Malchiodi}.

\begin{proposition}
\label{Taylor.expansion.of.scalar.curvature}
Consider a Riemannian metric on $\mathbb{R}^n$ of the form $g(x) =
\exp(h(x))$, where $h(x)$ is a trace-free symmetric two-tensor on
$\mathbb{R}^n$ satisfying $|h(x)| \leq 1$ for all $x \in \mathbb{R}^n$.
Let $R_g$ be the scalar curvature of $g$. There exists a constant $C$,
depending only on $n$, such that
\begin{align*}
&\Big | R_g - \partial_i \partial_k h_{ik}
+ \partial_i(h_{il} \, \partial_k h_{kl}) - \frac{1}{2} \, \partial_i
h_{il} \, \partial_k h_{kl} + \frac{1}{4} \, \partial_l h_{ik} \,
\partial_l h_{ik} \Big | \\
&\leq C \, |h|^2 \, |\partial^2 h| + C \, |h| \, |\partial h|^2.
\end{align*}
\end{proposition}

\vspace{2mm}

\begin{proposition}
\label{estimate.for.error.term}
Assume that $(\xi,\varepsilon) \in \lambda \, \Omega$. Then we have
\begin{align*}
&\Big \| \Delta_g u_{(\xi,\varepsilon)}
- \frac{n-2}{4(n-1)} \, R_g \, u_{(\xi,\varepsilon)}
+ n(n-2) \, u_{(\xi,\varepsilon)}^{\frac{n+2}{n-2}} \Big
\|_{L^{\frac{2n}{n+2}}(\mathbb{R}^n)} \\
&\leq C \, \lambda^8 \, \mu + C \, \Big ( \frac{\lambda}{\rho} \Big
)^{\frac{n-2}{2}}
\end{align*}
and
\begin{align*}
&\Big \| \Delta_g u_{(\xi,\varepsilon)}
- \frac{n-2}{4(n-1)} \, R_g \, u_{(\xi,\varepsilon)}
+ n(n-2) \, u_{(\xi,\varepsilon)}^{\frac{n+2}{n-2}} \\
&\hspace{10mm} + \sum_{i,k=1}^n \mu \, \lambda^6 \, f(\lambda^{-2} \,
|x|^2) \, H_{ik}(x) \, \partial_i \partial_k u_{(\xi,\varepsilon)} \Big
\|_{L^{\frac{2n}{n+2}}(\mathbb{R}^n)} \\
&\leq C \, \lambda^{\frac{8(n+2)}{n-2}} \, \mu^2 + C \, \Big (
\frac{\lambda}{\rho} \Big )^{\frac{n-2}{2}}.
\end{align*}
\end{proposition}

\textbf{Proof.} 
Note that $\sum_{i=1}^n \partial_i h_{ik}(x) = 0$ for $|x| \leq \rho$.
Hence, it follows from Proposition
\ref{Taylor.expansion.of.scalar.curvature} that
\[|R_g(x)| \leq C \, |h(x)|^2 \, |\partial^2 h(x)| + C \, |\partial
h(x)|^2 \leq C \, \mu^2 \, (\lambda + |x|)^{14}\]
for $|x| \leq \rho$. This implies
\begin{align*}
&\Big | \Delta_g u_{(\xi,\varepsilon)} - \frac{n-2}{4(n-1)} \, R_g \,
u_{(\xi,\varepsilon)}
+ n(n-2) \, u_{(\xi,\varepsilon)}^{\frac{n+2}{n-2}} \Big | \\
&= \Big | \sum_{i,k=1}^n \partial_i \big [ (g^{ik} - \delta_{ik}) \,
\partial_k u_{(\xi,\varepsilon)} \big ]
- \frac{n-2}{4(n-1)} \, R_g \, u_{(\xi,\varepsilon)} \Big | \\
&\leq C \, \lambda^{\frac{n-2}{2}} \, \mu \, (\lambda + |x|)^{8-n}
\end{align*}
and
\begin{align*}
&\Big | \Delta_g u_{(\xi,\varepsilon)}
- \frac{n-2}{4(n-1)} \, R_g \, u_{(\xi,\varepsilon)}
+ n(n-2) \, u_{(\xi,\varepsilon)}^{\frac{n+2}{n-2}} + \sum_{i,k=1}^n
h_{ik} \, \partial_i \partial_k u_{(\xi,\varepsilon)} \Big | \\
&= \Big | \sum_{i,k=1}^n \partial_i \big [ (g^{ik} - \delta_{ik} +
h_{ik}) \, \partial_k u_{(\xi,\varepsilon)} \big ] - \frac{n-2}{4(n-1)}
\, R_g \, u_{(\xi,\varepsilon)} \Big | \\
&\leq C \, \lambda^{\frac{n-2}{2}} \, \mu^2 \, (\lambda + |x|)^{16-n} \\
&\leq C \, \lambda^{\frac{n-2}{2}} \, \mu^2 \, (\lambda +
|x|)^{\frac{8(n+2)}{n-2}-n}
\end{align*}
for $|x| \leq \rho$. From this the assertion follows. \\

\begin{corollary}
\label{estimate.for.v.1}
The function $v_{(\xi,\varepsilon)} - u_{(\xi,\varepsilon)}$ satisfies
the estimate
\[\|v_{(\xi,\varepsilon)} -
u_{(\xi,\varepsilon)}\|_{L^{\frac{2n}{n-2}}(\mathbb{R}^n)}
\leq C \, \lambda^8 \, \mu + C \, \Big ( \frac{\lambda}{\rho} \Big
)^{\frac{n-2}{2}}\]
whenever $(\xi,\varepsilon) \in \lambda \, \Omega$.
\end{corollary}

\textbf{Proof.}
It follows from Proposition \ref{fixed.point.argument} that
\begin{align*}
&\|v_{(\xi,\varepsilon)} -
u_{(\xi,\varepsilon)}\|_{L^{\frac{2n}{n-2}}(\mathbb{R}^n)} \\
&\leq C \, \Big \| \Delta_g u_{(\xi,\varepsilon)}
- \frac{n-2}{4(n-1)} \, R_g \, u_{(\xi,\varepsilon)}
+ n(n-2) \, u_{(\xi,\varepsilon)}^{\frac{n+2}{n-2}} \Big
\|_{L^{\frac{2n}{n+2}}(\mathbb{R}^n)},
\end{align*}
where $C$ is a constant that depends only on $n$. Hence, the assertion
follows from Proposition \ref{estimate.for.error.term}. \\

We next establish a more precise estimate for the function $v_{(\xi,\varepsilon)} - u_{(\xi,\varepsilon)}$. 
Applying Proposition \ref{linearized.operator} with $h = 0$, we conclude that there exists a unique function $w_{(\xi,\varepsilon)} \in
\mathcal{E}_{(\xi,\varepsilon)}$ such that
\begin{align*}
&\int_{\mathbb{R}^n} \Big ( \langle dw_{(\xi,\varepsilon)},d\psi \rangle
- n(n+2) \, u_{(\xi,\varepsilon)}^{\frac{4}{n-2}} \,
w_{(\xi,\varepsilon)} \, \psi \Big ) \\
&= -\int_{\mathbb{R}^n} \sum_{i,k=1}^n \mu \, \lambda^6 \,
f(\lambda^{-2} \, |x|^2) \, H_{ik}(x) \, \partial_i \partial_k
u_{(\xi,\varepsilon)} \, \psi
\end{align*}
for all test functions $\psi \in \mathcal{E}_{(\xi,\varepsilon)}$. 

\begin{proposition} 
\label{properties.of.w}
The function $w_{(\xi,\varepsilon)}$ is smooth. Moreover, if $(\xi,\varepsilon) \in \lambda \, \Omega$, then the function
$w_{(\xi,\varepsilon)}$ satisfies the estimates
\begin{align*}
&|w_{(\xi,\varepsilon)}(x)| \leq C \, \lambda^{\frac{n-2}{2}} \, \mu \,
(\lambda + |x|)^{10-n} \\
&|\partial w_{(\xi,\varepsilon)}(x)| \leq C \, \lambda^{\frac{n-2}{2}}
\, \mu \, (\lambda + |x|)^{9-n} \\
&|\partial^2 w_{(\xi,\varepsilon)}(x)| \leq C \, \lambda^{\frac{n-2}{2}}
\, \mu \, (\lambda + |x|)^{8-n}
\end{align*}
for all $x \in \mathbb{R}^n$.
\end{proposition}

\textbf{Proof.} 
There exist real numbers $b_k(\xi,\varepsilon)$ such that 
\begin{align*}
&\int_{\mathbb{R}^n} \Big ( \langle dw_{(\xi,\varepsilon)},d\psi \rangle
- n(n+2) \, u_{(\xi,\varepsilon)}^{\frac{4}{n-2}} \,
w_{(\xi,\varepsilon)} \, \psi \Big ) \\
&= -\int_{\mathbb{R}^n} \sum_{i,k=1}^n \mu \, \lambda^6 \,
f(\lambda^{-2} \, |x|^2) \, H_{ik}(x) \, \partial_i \partial_k
u_{(\xi,\varepsilon)} \, \psi \\ 
&+ \sum_{k=0}^n b_k(\xi,\varepsilon) \, \int_{\mathbb{R}^n} \varphi_{(\xi,\varepsilon,k)} \, \psi
\end{align*}
for all test functions $\psi \in \mathcal{E}$. Hence, standard elliptic regularity theory implies that $w_{(\xi,\varepsilon)}$ is smooth. 

It remains to prove quantitative estimates for $w_{(\xi,\varepsilon)}$. To that end, we consider a pair $(\xi,\varepsilon) \in \lambda \, \Omega$. One readily verifies that 
\[\Big \| \sum_{i,k=1}^n \mu \, \lambda^6 \,
f(\lambda^{-2} \, |x|^2) \, H_{ik}(x) \, \partial_i \partial_k
u_{(\xi,\varepsilon)} \Big \|_{L^{\frac{2n}{n+2}}(\mathbb{R}^n)} \leq C \, \lambda^8 \, \mu.\] As a consequence, the function $w_{(\xi,\varepsilon)}$ satisfies $\|w_{(\xi,\varepsilon)}\|_{L^{\frac{2n}{n-2}}(\mathbb{R}^n)} \leq C \, \lambda^8 \, \mu$. Moreover, we have $\sum_{k=0}^n |b_k(\xi,\varepsilon)| \leq C \, \lambda^8 \, \mu$. This implies 
\begin{align*} 
&\big | \Delta w_{(\xi,\varepsilon)} + n(n+2) \, u_{(\xi,\varepsilon)}^{\frac{4}{n-2}} \, w_{(\xi,\varepsilon)} \big | \\ 
&= \bigg | \sum_{i,k=1}^n \mu \, \lambda^6 \,
f(\lambda^{-2} \, |x|^2) \, H_{ik}(x) \, \partial_i \partial_k
u_{(\xi,\varepsilon)} - \sum_{k=0}^n b_k(\xi,\varepsilon) \, \int_{\mathbb{R}^n} \varphi_{(\xi,\varepsilon,k)} \bigg | \\ 
&\leq C \, \lambda^{\frac{n-2}{2}} \, \mu \, (\lambda+|x|)^{8-n} 
\end{align*}
for all $x \in \mathbb{R}^n$. We claim that 
\[\sup_{x \in \mathbb{R}^n} (\lambda + |x|)^{\frac{n-2}{2}} \, |w_{(\xi,\varepsilon)}(x)| \leq C \, \lambda^8 \, \mu.\] 
To show this, we fix a point $x_0 \in \mathbb{R}^n$. Let $r = \frac{1}{2} \, (\lambda + |x_0|)$. Then 
\[u_{(\xi,\varepsilon)}(x)^{\frac{4}{n-2}} \leq C \, r^{-2}\] and 
\[\big | \Delta w_{(\xi,\varepsilon)} + n(n+2) \, u_{(\xi,\varepsilon)}^{\frac{4}{n-2}} \, w_{(\xi,\varepsilon)} \big | \leq C \, \lambda^{\frac{n-2}{2}} \, \mu \, r^{8-n}\] 
for all $x \in B_r(x_0)$. Hence, it follows from standard interior estimates that 
\begin{align*} 
r^{\frac{n-2}{2}} \, |w_{(\xi,\varepsilon)}(x_0)| &\leq C \, \|w_{(\xi,\varepsilon)}\|_{L^{\frac{2n}{n-2}}(B_r(x_0))} \\ 
&+ C \, r^{\frac{n+2}{2}} \, \big \| \Delta w_{(\xi,\varepsilon)} + n(n+2) \, u_{(\xi,\varepsilon)}^{\frac{4}{n-2}} \, w_{(\xi,\varepsilon)} \big \|_{L^\infty(B_r(x_0))} \\ 
&\leq C \, \lambda^8 \, \mu + C \, \lambda^{\frac{n-2}{2}} \, \mu \, r^{-\frac{n-18}{2}} \\ 
&\leq C \, \lambda^8 \, \mu. 
\end{align*} 
Therefore, we have 
\[\sup_{x \in \mathbb{R}^n} (\lambda + |x|)^{\frac{n-2}{2}} \, |w_{(\xi,\varepsilon)}(x)| \leq C \, \lambda^8 \, \mu,\] as 
claimed. Since $\sup_{x \in \mathbb{R}^n} |x|^{\frac{n-2}{2}} \, |w_{(\xi,\varepsilon)}(x)| < \infty$, we can express the function $w_{(\xi,\varepsilon)}$ in the form 
\begin{equation} 
\label{convolution.formula}
w_{(\xi,\varepsilon)}(x) = -\frac{1}{(n-2) \, |S^{n-1}|} \int_{\mathbb{R}^n} |x - y|^{2-n} \, \Delta w_{(\xi,\varepsilon)}(y) \, dy 
\end{equation}
for all $x \in \mathbb{R}^n$. 

We are now able to use a bootstrap argument to prove the desired estimate for $w_{(\xi,\varepsilon)}$. It follows from (\ref{convolution.formula}) that
\[\sup_{x \in \mathbb{R}^n} (\lambda + |x|)^\beta \, |w_{(\xi,\varepsilon)}(x)| \leq C \, \sup_{x \in \mathbb{R}^n} (\lambda + |x|)^{\beta+2} \, |\Delta w_{(\xi,\varepsilon)}(x)|\]
for all $0 < \beta < n-2$. Since 
\begin{align*} 
|\Delta w_{(\xi,\varepsilon)}(x)| 
&\leq n(n+2) \, u_{(\xi,\varepsilon)}(x)^{\frac{4}{n-2}} \, |w_{(\xi,\varepsilon)}(x)| \\ 
&+ C \, \lambda^{\frac{n-2}{2}} \, \mu \, (\lambda+|x|)^{8-n} 
\end{align*}
for all $x \in \mathbb{R}^n$, we conclude that
\begin{align*}
\sup_{x \in \mathbb{R}^n} (\lambda + |x|)^\beta \, |w_{(\xi,\varepsilon)}(x)| &\leq C \, \lambda^2 \, \sup_{x \in \mathbb{R}^n} (\lambda + |x|)^{\beta-2} \, |w_{(\xi,\varepsilon)}(x)| \\ &+ C \, \lambda^{\beta-\frac{n-18}{2}} \, \mu
\end{align*}
for all $0 < \beta \leq n-10$. Iterating this inequality, we obtain
\[\sup_{x \in \mathbb{R}^n} (\lambda + |x|)^{n-10} \, |w_{(\xi,\varepsilon)}(x)| \leq C \, \lambda^{\frac{n-2}{2}} \, \mu.\] 
The estimates for the first and second derivatives of $w_{(\xi,\varepsilon)}$ follow now from standard interior estimates. \\

\begin{corollary}
\label{estimate.for.v.2}
The function $v_{(\xi,\varepsilon)} - u_{(\xi,\varepsilon)} -
w_{(\xi,\varepsilon)}$ satisfies the estimate
\[\|v_{(\xi,\varepsilon)} - u_{(\xi,\varepsilon)} -
w_{(\xi,\varepsilon)}\|_{L^{\frac{2n}{n-2}}(\mathbb{R}^n)}
\leq C \, \lambda^{\frac{8(n+2)}{n-2}} \, \mu^{\frac{n+2}{n-2}} + C \,
\Big ( \frac{\lambda}{\rho} \Big )^{\frac{n-2}{2}}\]
whenever $(\xi,\varepsilon) \in \lambda \, \Omega$.
\end{corollary}

\textbf{Proof.}
Consider the functions
\[B_1 = \sum_{i,k=1}^n \partial_i \big [ (g^{ik} - \delta_{ik}) \,
\partial_k w_{(\xi,\varepsilon)} \big ] - \frac{n-2}{4(n-1)} \, R_g \,
w_{(\xi,\varepsilon)}\]
and
\[B_2 = \sum_{i,k=1}^n \mu \, \lambda^6 \, f(\lambda^{-2} \, |x|^2) \,
H_{ik}(x) \, \partial_i \partial_k u_{(\xi,\varepsilon)}.\]
By definition of $w_{(\xi,\varepsilon)}$, we have
\begin{align*}
&\int_{\mathbb{R}^n} \Big ( \langle dw_{(\xi,\varepsilon)},d\psi \rangle_g
+ \frac{n-2}{4(n-1)} \, R_g \, w_{(\xi,\varepsilon)} \, \psi - n(n+2) \,
u_{(\xi,\varepsilon)}^{\frac{4}{n-2}} \, w_{(\xi,\varepsilon)} \, \psi
\Big ) \\
&= -\int_{\mathbb{R}^n} (B_1 + B_2) \, \psi
\end{align*}
for all functions $\psi \in \mathcal{E}_{(\xi,\varepsilon)}$. Since
$w_{(\xi,\varepsilon)} \in \mathcal{E}_{(\xi,\varepsilon)}$, it follows
that
\[w_{(\xi,\varepsilon)} = -G_{(\xi,\varepsilon)}(B_1 + B_2).\]
Moreover, we have
\[v_{(\xi,\varepsilon)} - u_{(\xi,\varepsilon)} = G_{(\xi,\varepsilon)}
\big ( B_3 + n(n-2) \, B_4 \big ),\]
where
\[B_3 = \Delta_g u_{(\xi,\varepsilon)} - \frac{n-2}{4(n-1)} \, R_g \,
u_{(\xi,\varepsilon)} + n(n-2) \, u_{(\xi,\varepsilon)}^{\frac{n+2}{n-2}}\]
and
\[B_4 = |v_{(\xi,\varepsilon)}|^{\frac{4}{n-2}} \, v_{(\xi,\varepsilon)}
- u_{(\xi,\varepsilon)}^{\frac{n+2}{n-2}}  - \frac{n+2}{n-2} \,
u_{(\xi,\varepsilon)}^{\frac{4}{n-2}} \, (v_{(\xi,\varepsilon)} -
u_{(\xi,\varepsilon)}).\]
Thus, we conclude that
\[v_{(\xi,\varepsilon)} - u_{(\xi,\varepsilon)} - w_{(\xi,\varepsilon)}
= G_{(\xi,\varepsilon)} \big ( B_1 + B_2 + B_3 + n(n-2) \, B_4 \big ),\]
where $G_{(\xi,\varepsilon)}$ denotes the solution operator constructed in Proposition \ref{linearized.operator}. As a consequence, we obtain 
\[\|v_{(\xi,\varepsilon)} - u_{(\xi,\varepsilon)} - w_{(\xi,\varepsilon)}\|_{L^{\frac{2n}{n-2}}(\mathbb{R}^n)} \leq C \, 
\big \| B_1 + B_2 + B_3 + n(n-2) \, B_4 \big \|_{L^{\frac{2n}{n+2}}(\mathbb{R}^n)}.\]
It follows from Proposition \ref{properties.of.w} that 
\[|B_1(x)| \leq C \, \lambda^{\frac{n-2}{2}} \, \mu^2 \, (\lambda + |x|)^{16-n} \leq C \, \lambda^{\frac{n-2}{2}} \, \mu^2 \, (\lambda+|x|)^{\frac{8(n+2)}{n-2}-n}\] 
for $|x| \leq \rho$ and 
\[|B_1(x)| \leq C \, \lambda^{\frac{n-2}{2}} \, \mu \, |x|^{8-n}\] 
for $|x| \geq \rho$. This implies 
\[\|B_1\|_{L^{\frac{2n}{n+2}}(\mathbb{R}^n)} \leq C \,
\lambda^{\frac{8(n+2)}{n-2}} \, \mu^2 + C \, \rho^8 \, \mu \, \Big (
\frac{\lambda}{\rho} \Big )^{\frac{n-2}{2}}.\]
Moreover, observe that 
\[\|B_2 + B_3\|_{L^{\frac{2n}{n+2}}(\mathbb{R}^n)} \leq C \,
\lambda^{\frac{8(n+2)}{n-2}} \, \mu^2 + C \, \Big ( \frac{\lambda}{\rho}
\Big )^{\frac{n-2}{2}}\]
by Proposition \ref{estimate.for.error.term}. Finally, Corollary \ref{estimate.for.v.1} implies that 
\begin{align*}
\|B_4\|_{L^{\frac{2n}{n+2}}(\mathbb{R}^n)}
&\leq C \, \|v_{(\xi,\varepsilon)} -
u_{(\xi,\varepsilon)}\|_{L^{\frac{2n}{n-2}}(\mathbb{R}^n)}^{\frac{n+2}{n-2}}
\\
&\leq C \, \lambda^{\frac{8(n+2)}{n-2}} \, \mu^{\frac{n+2}{n-2}} + C \,
\Big ( \frac{\lambda}{\rho} \Big )^{\frac{n+2}{2}}.
\end{align*}
Putting these facts together, we obtain 
\[\|v_{(\xi,\varepsilon)} - u_{(\xi,\varepsilon)} -
w_{(\xi,\varepsilon)}\|_{L^{\frac{2n}{n-2}}(\mathbb{R}^n)}
\leq C \, \lambda^{\frac{8(n+2)}{n-2}} \, \mu^{\frac{n+2}{n-2}} + C \,
\Big ( \frac{\lambda}{\rho} \Big )^{\frac{n-2}{2}}.
\]
This completes the proof. \\

\begin{proposition}
\label{term.1}
We have
\begin{align*}
&\bigg | \int_{\mathbb{R}^n} \Big ( |dv_{(\xi,\varepsilon)}|_g^2
- |du_{(\xi,\varepsilon)}|_g^2 + \frac{n-2}{4(n-1)} \, R_g \,
(v_{(\xi,\varepsilon)}^2 - u_{(\xi,\varepsilon)}^2) \Big ) \\
&\hspace{10mm} + \int_{\mathbb{R}^n} n(n-2) \,
(|v_{(\xi,\varepsilon)}|^{\frac{4}{n-2}} -
u_{(\xi,\varepsilon)}^{\frac{4}{n-2}}) \, u_{(\xi,\varepsilon)} \,
v_{(\xi,\varepsilon)} \\
&\hspace{10mm} - \int_{\mathbb{R}^n} n(n-2) \,
(|v_{(\xi,\varepsilon)}|^{\frac{2n}{n-2}} -
u_{(\xi,\varepsilon)}^{\frac{2n}{n-2}}) \\
&\hspace{10mm} - \int_{\mathbb{R}^n} \sum_{i,k=1}^n \mu \, \lambda^6 \,
f(\lambda^{-2} \, |x|^2) \, H_{ik}(x) \, \partial_i \partial_k
u_{(\xi,\varepsilon)} \, w_{(\xi,\varepsilon)} \bigg | \\
&\leq C \, \lambda^{\frac{16n}{n-2}} \, \mu^{\frac{2n}{n-2}} + C \,
\lambda^8 \, \mu \, \Big ( \frac{\lambda}{\rho} \Big )^{\frac{n-2}{2}} +
C \, \Big ( \frac{\lambda}{\rho} \Big )^{n-2}
\end{align*}
whenever $(\xi,\varepsilon) \in \lambda \, \Omega$.
\end{proposition}

\textbf{Proof.}
By definition of $v_{(\xi,\varepsilon)}$, we have
\begin{align*}
&\int_{\mathbb{R}^n} \Big ( |dv_{(\xi,\varepsilon)}|_g^2 - \langle
du_{(\xi,\varepsilon)},dv_{(\xi,\varepsilon)} \rangle_g
+ \frac{n-2}{4(n-1)} \, R_g \, v_{(\xi,\varepsilon)} \,
(v_{(\xi,\varepsilon)} - u_{(\xi,\varepsilon)}) \Big ) \\
&\hspace{10mm} - \int_{\mathbb{R}^n} n(n-2) \,
|v_{(\xi,\varepsilon)}|^{\frac{4}{n-2}} \, v_{(\xi,\varepsilon)} \,
(v_{(\xi,\varepsilon)} - u_{(\xi,\varepsilon)}) = 0.
\end{align*}
Using Proposition \ref{estimate.for.error.term} and Corollary
\ref{estimate.for.v.1}, we obtain
\begin{align*}
&\bigg | \int_{\mathbb{R}^n} \Big ( \langle
du_{(\xi,\varepsilon)},dv_{(\xi,\varepsilon)} \rangle_g
- |du_{(\xi,\varepsilon)}|_g^2 + \frac{n-2}{4(n-1)} \, R_g \,
u_{(\xi,\varepsilon)} \, (v_{(\xi,\varepsilon)} - u_{(\xi,\varepsilon)})
\Big ) \\
&\hspace{10mm} - \int_{\mathbb{R}^n} n(n-2) \,
u_{(\xi,\varepsilon)}^{\frac{n+2}{n-2}} \, (v_{(\xi,\varepsilon)} -
u_{(\xi,\varepsilon)}) \\
&\hspace{10mm} - \int_{\mathbb{R}^n} \sum_{i,k=1}^n \mu \, \lambda^6 \,
f(\lambda^{-2} \, |x|^2) \, H_{ik}(x) \, \partial_i \partial_k
u_{(\xi,\varepsilon)} \, (v_{(\xi,\varepsilon)} - u_{(\xi,\varepsilon)})
\bigg | \\
&\leq \Big \| \Delta_g u_{(\xi,\varepsilon)}
- \frac{n-2}{4(n-1)} \, R_g \, u_{(\xi,\varepsilon)}
+ n(n-2) \, u_{(\xi,\varepsilon)}^{\frac{n+2}{n-2}} \\
&\hspace{10mm} + \sum_{i,k=1}^n \mu \, \lambda^6 \, f(\lambda^{-2} \,
|x|^2) \, H_{ik}(x) \, \partial_i \partial_k u_{(\xi,\varepsilon)} \Big
\|_{L^{\frac{2n}{n+2}}(\mathbb{R}^n)} \\
&\hspace{5mm} \cdot \|v_{(\xi,\varepsilon)} -
u_{(\xi,\varepsilon)}\|_{L^{\frac{2n}{n-2}}(\mathbb{R}^n)} \\
&\leq C \, \lambda^{\frac{16n}{n-2}} \, \mu^3 + C \, \lambda^8 \, \mu \,
\Big ( \frac{\lambda}{\rho} \Big )^{\frac{n-2}{2}} + C \, \Big (
\frac{\lambda}{\rho} \Big )^{n-2}.
\end{align*}
Moreover, we have
\begin{align*}
&\bigg | \int_{\mathbb{R}^n} \sum_{i,k=1}^n \mu \, \lambda^6 \,
f(\lambda^{-2} \, |x|^2) \, H_{ik}(x) \, \partial_i \partial_k
u_{(\xi,\varepsilon)} \, (v_{(\xi,\varepsilon)} - u_{(\xi,\varepsilon)}
- w_{(\xi,\varepsilon)}) \bigg | \\
&\leq C \, \lambda^8 \, \mu \, \|v_{(\xi,\varepsilon)} -
u_{(\xi,\varepsilon)} -
w_{(\xi,\varepsilon)}\|_{L^{\frac{2n}{n-2}}(\mathbb{R}^n)} \\
&\leq C \, \lambda^{\frac{16n}{n-2}} \, \mu^{\frac{2n}{n-2}} + C \,
\lambda^8 \, \mu \, \Big ( \frac{\lambda}{\rho} \Big )^{\frac{n-2}{2}}
\end{align*}
by Corollary \ref{estimate.for.v.2}. Putting these facts together, the
assertion follows. \\

\begin{proposition}
\label{term.2}
We have
\begin{align*}
&\bigg | \int_{\mathbb{R}^n} (|v_{(\xi,\varepsilon)}|^{\frac{4}{n-2}} -
u_{(\xi,\varepsilon)}^{\frac{4}{n-2}}) \, u_{(\xi,\varepsilon)} \,
v_{(\xi,\varepsilon)} - \frac{2}{n} \int_{\mathbb{R}^n}
(|v_{(\xi,\varepsilon)}|^{\frac{2n}{n-2}} -
u_{(\xi,\varepsilon)}^{\frac{2n}{n-2}}) \bigg | \\
&\leq C \, \lambda^{\frac{16n}{n-2}} \, \mu^{\frac{2n}{n-2}} + C \, \Big
( \frac{\lambda}{\rho} \Big )^n
\end{align*}
whenever $(\xi,\varepsilon) \in \lambda \, \Omega$.
\end{proposition}

\textbf{Proof.}
We have the pointwise estimate
\begin{align*}
&\Big | (|v_{(\xi,\varepsilon)}|^{\frac{4}{n-2}} -
u_{(\xi,\varepsilon)}^{\frac{4}{n-2}}) \, u_{(\xi,\varepsilon)} \,
v_{(\xi,\varepsilon)} - \frac{2}{n} \,
(|v_{(\xi,\varepsilon)}|^{\frac{2n}{n-2}} -
u_{(\xi,\varepsilon)}^{\frac{2n}{n-2}}) \Big | \\
&\leq C \, |v_{(\xi,\varepsilon)} -
u_{(\xi,\varepsilon)}|^{\frac{2n}{n-2}},
\end{align*}
where $C$ is a constant that depends only on $n$. This implies
\begin{align*}
&\bigg | \int_{\mathbb{R}^n} (|v_{(\xi,\varepsilon)}|^{\frac{4}{n-2}} -
u_{(\xi,\varepsilon)}^{\frac{4}{n-2}}) \, u_{(\xi,\varepsilon)} \,
v_{(\xi,\varepsilon)} - \frac{2}{n} \int_{\mathbb{R}^n}
(|v_{(\xi,\varepsilon)}|^{\frac{2n}{n-2}} -
u_{(\xi,\varepsilon)}^{\frac{2n}{n-2}}) \bigg | \\
&\leq C \, \|v_{(\xi,\varepsilon)} -
u_{(\xi,\varepsilon)}\|_{L^{\frac{2n}{n-2}}(\mathbb{R}^n)}^{\frac{2n}{n-2}}
\\
&\leq C \, \lambda^{\frac{16n}{n-2}} \, \mu^{\frac{2n}{n-2}} + C \, \Big
( \frac{\lambda}{\rho} \Big )^n 
\end{align*}
by Corollary \ref{estimate.for.v.1}. \\

\begin{proposition}
\label{term.3}
We have
\begin{align*}
&\bigg | \int_{\mathbb{R}^n} \Big ( |du_{(\xi,\varepsilon)}|_g^2 +
\frac{n-2}{4(n-1)} \, R_g \, u_{(\xi,\varepsilon)}^2 - n(n-2) \,
u_{(\xi,\varepsilon)}^{\frac{2n}{n-2}} \Big ) \\
&\hspace{10mm} - \int_{B_\rho(0)} \frac{1}{2} \, \sum_{i,k,l=1}^n
h_{il} \, h_{kl} \, \partial_i u_{(\xi,\varepsilon)} \, \partial_k
u_{(\xi,\varepsilon)} \\
&\hspace{10mm} + \int_{B_\rho(0)} \frac{n-2}{16(n-1)} \,
\sum_{i,k,l=1}^n (\partial_l h_{ik})^2 \,
u_{(\xi,\varepsilon)}^2 \bigg | \\
&\leq C \, \lambda^{\frac{16n}{n-2}} \, \mu^3 + C \, \Big (
\frac{\lambda}{\rho} \Big )^{n-2}
\end{align*}
whenever $(\xi,\varepsilon) \in \lambda \, \Omega$.
\end{proposition}

\textbf{Proof.}
Note that
\begin{align*}
&\Big | g^{ik}(x) - \delta_{ik} + h_{ik}(x) - \frac{1}{2} \,
\sum_{l=1}^n h_{il}(x) \, h_{kl}(x) \Big | \\
&\leq C \, |h(x)|^3 \leq C \, \mu^3 \, (\lambda + |x|)^{24} \leq C \,
\mu^3 \, (\lambda+|x|)^{\frac{16n}{n-2}}
\end{align*}
for $|x| \leq \rho$. This implies
\begin{align*}
&\bigg | \int_{\mathbb{R}^n} \big ( |du_{(\xi,\varepsilon)}|_g^2 -
|du_{(\xi,\varepsilon)}|^2 \big ) + \int_{\mathbb{R}^n} \sum_{i,k=1}^n
h_{ik} \, \partial_i u_{(\xi,\varepsilon)} \, \partial_k u_{(\xi,\varepsilon)} \\
&\hspace{10mm} - \int_{B_\rho(0)} \frac{1}{2} \sum_{i,k,l=1}^n h_{il}
\, h_{kl} \, \partial_i u_{(\xi,\varepsilon)} \, \partial_k
u_{(\xi,\varepsilon)} \bigg | \\
&\leq C \, \lambda^{n-2} \, \mu^3 \, \int_{B_\rho(0)} (\lambda+|x|)^{\frac{16n}{n-2}+2-2n} 
+ C \, \lambda^{n-2} \, \int_{\mathbb{R}^n \setminus B_\rho(0)} (\lambda + |x|)^{2-2n} \\ 
&\leq C \, \lambda^{\frac{16n}{n-2}} \, \mu^3 + C \, \Big ( \frac{\lambda}{\rho} \Big )^{n-2}.
\end{align*}
By Proposition \ref{Taylor.expansion.of.scalar.curvature}, the scalar
curvature of $g$ satisfies the estimate
\begin{align*}
&\Big | R_g(x) + \frac{1}{4} \sum_{i,k,l=1}^n (\partial_l h_{ik}(x))^2
\Big | \\
&\leq C \, |h(x)|^2 \, |\partial^2 h(x)| + C \, |h(x)| \, |\partial
h(x)|^2 \\
&\leq C \, \mu^3 \, (\lambda + |x|)^{22} \leq C \, \mu^3 \, (\lambda +
|x|)^{\frac{16n}{n-2}-2}
\end{align*}
for $|x| \leq \rho$. This implies
\begin{align*}
&\bigg | \int_{\mathbb{R}^n} R_g \, u_{(\xi,\varepsilon)}^2 +
\int_{B_\rho(0)} \frac{1}{4} \sum_{i,k,l=1}^n (\partial_l h_{ik})^2
\, u_{(\xi,\varepsilon)}^2 \bigg | \\
&\leq C \, \lambda^{n-2} \, \mu^3 \, \int_{B_\rho(0)} (\lambda+|x|)^{\frac{16n}{n-2}+2-2n} 
+ C \, \lambda^{n-2} \, \int_{\mathbb{R}^n \setminus B_\rho(0)} (\lambda + |x|)^{4-2n} \\ 
&\leq C \, \lambda^{\frac{16n}{n-2}} \, \mu^3 + C \, \rho^2 \, \Big (
\frac{\lambda}{\rho} \Big )^{n-2}.
\end{align*} 
At this point, we use the formula
\begin{align*}
&\partial_i u_{(\xi,\varepsilon)} \, \partial_k u_{(\xi,\varepsilon)} -
\frac{n-2}{4(n-1)} \, \partial_i \partial_k (u_{(\xi,\varepsilon)}^2) \\
&= \frac{1}{n} \, \Big ( |du_{(\xi,\varepsilon)}|^2 - \frac{n-2}{4(n-1)}
\, \Delta (u_{(\xi,\varepsilon)}^2) \Big ) \, \delta_{ik}.
\end{align*}
Since $h_{ik}$ is trace-free, we obtain
\[\sum_{i,k=1}^n h_{ik} \, \partial_i u_{(\xi,\varepsilon)} \,
\partial_k u_{(\xi,\varepsilon)} = \frac{n-2}{4(n-1)} \sum_{i,k=1}^n
h_{ik} \, \partial_i \partial_k (u_{(\xi,\varepsilon)}^2),\]
hence
\[\int_{\mathbb{R}^n} \sum_{i,k=1}^n h_{ik} \, \partial_i
u_{(\xi,\varepsilon)} \, \partial_k u_{(\xi,\varepsilon)} =
\int_{\mathbb{R}^n} \frac{n-2}{4(n-1)} \sum_{i,k=1}^n \partial_i
\partial_k h_{ik} \, u_{(\xi,\varepsilon)}^2.\]
Since $\sum_{i=1}^n \partial_i h_{ik}(x) = 0$ for $|x| \leq \rho$, it
follows that
\[\bigg | \int_{\mathbb{R}^n} \sum_{i,k=1}^n h_{ik} \, \partial_i
u_{(\xi,\varepsilon)} \, \partial_k u_{(\xi,\varepsilon)} \bigg |
\leq C \int_{\mathbb{R}^n \setminus B_\rho(0)} u_{(\xi,\varepsilon)}^2 \leq C \, \rho^2 \, \Big ( \frac{\lambda}{\rho} \Big )^{n-2}.\]
Putting these facts together, the assertion follows. \\

\begin{corollary}
\label{key.estimate}
The function $\mathcal{F}_g(\xi,\varepsilon)$ satisfies the estimate
\begin{align*}
&\bigg | \mathcal{F}_g(\xi,\varepsilon) - \int_{B_\rho(0)}
\frac{1}{2} \, \sum_{i,k,l=1}^n h_{il} \, h_{kl} \, \partial_i
u_{(\xi,\varepsilon)} \, \partial_k u_{(\xi,\varepsilon)} \\
&\hspace{10mm} + \int_{B_\rho(0)} \frac{n-2}{16(n-1)} \,
\sum_{i,k,l=1}^n (\partial_l h_{ik})^2 \,
u_{(\xi,\varepsilon)}^2 \\
&\hspace{10mm} - \int_{\mathbb{R}^n} \sum_{i,k=1}^n \mu \, \lambda^6 \,
f(\lambda^{-2} \, |x|^2) \, H_{ik}(x) \, \partial_i\partial_k
u_{(\xi,\varepsilon)} \, w_{(\xi,\varepsilon)} \bigg | \\
&\leq C \, \lambda^{\frac{16n}{n-2}} \, \mu^{\frac{2n}{n-2}}
+ C \, \lambda^8 \, \mu \, \Big ( \frac{\lambda}{\rho} \Big
)^{\frac{n-2}{2}} + C \, \Big ( \frac{\lambda}{\rho} \Big )^{n-2}
\end{align*}
whenever $(\xi,\varepsilon) \in \lambda \, \Omega$.
\end{corollary}

\section{Finding a critical point of an auxiliary function}

We define a function $F: \mathbb{R}^n \times (0,\infty) \to \mathbb{R}$
as follows: given any pair $(\xi,\varepsilon) \in \mathbb{R}^n \times
(0,\infty)$, we define
\begin{align*}
F(\xi,\varepsilon)
&= \int_{\mathbb{R}^n} \frac{1}{2} \sum_{i,k,l=1}^n \overline{H}_{il}(x)
\, \overline{H}_{kl}(x) \, \partial_i u_{(\xi,\varepsilon)}(x) \,
\partial_k u_{(\xi,\varepsilon)}(x) \\
&- \int_{\mathbb{R}^n} \frac{n-2}{16(n-1)} \, \sum_{i,k,l=1}^n
(\partial_l \overline{H}_{ik}(x))^2 \, u_{(\xi,\varepsilon)}(x)^2 \\
&+ \int_{\mathbb{R}^n} \sum_{i,k=1}^n \overline{H}_{ik}(x) \, \partial_i
\partial_k u_{(\xi,\varepsilon)}(x) \, z_{(\xi,\varepsilon)}(x),
\end{align*}
where $z_{(\xi,\varepsilon)} \in \mathcal{E}_{(\xi,\varepsilon)}$
satisfies the relation
\begin{align*}
&\int_{\mathbb{R}^n} \Big ( \langle dz_{(\xi,\varepsilon)},d\psi \rangle
- n(n+2) \, u_{(\xi,\varepsilon)}(x)^{\frac{4}{n-2}} \,
z_{(\xi,\varepsilon)} \, \psi \Big ) \\
&= -\int_{\mathbb{R}^n} \sum_{i,k=1}^n \overline{H}_{ik} \, \partial_i
\partial_k u_{(\xi,\varepsilon)} \, \psi
\end{align*}
for all test functions $\psi \in \mathcal{E}_{(\xi,\varepsilon)}$. Our
goal in this section is to show that the function $F(\xi,\varepsilon)$
has a critical point.

\begin{proposition}
\label{symmetry}
The function $F(\xi,\varepsilon)$ satisfies $F(\xi,\varepsilon) =
F(-\xi,\varepsilon)$ for all $(\xi,\varepsilon) \in \mathbb{R}^n \times
(0,\infty)$. Consequently, we have $\frac{\partial}{\partial \xi_p}
F(0,\varepsilon) = 0$ and $\frac{\partial^2}{\partial \varepsilon \,
\partial \xi_p} F(0,\varepsilon) = 0$ for all $\varepsilon > 0$ and $p =
1, \hdots, n$.
\end{proposition}

\textbf{Proof.}
This follows immediately from the relation $\overline{H}_{ik}(-x) =
\overline{H}_{ik}(x)$. \\

\begin{proposition}
\label{integral.identity.1}
We have
\begin{align*}
&\int_{\partial B_r(0)} \sum_{i,k,l=1}^n (\partial_l H_{ik}(x))^2 \, x_p
\, x_q \\
&= \frac{2}{n(n+2)} \, |S^{n-1}| \, \sum_{i,k,l=1}^n (W_{ipkl} +
W_{ilkp}) \, (W_{iqkl} + W_{ilkq}) \, r^{n+3} \\
&+ \frac{1}{n(n+2)} \, |S^{n-1}| \, \sum_{i,j,k,l=1}^n (W_{ijkl} +
W_{ilkj})^2 \, \delta_{pq} \, r^{n+3}
\end{align*}
and
\begin{align*}
&\int_{\partial B_r(0)} \sum_{i,k=1}^n H_{ik}(x)^2 \, x_p \, x_q \\
&= \frac{2}{n(n+2)(n+4)} \, |S^{n-1}| \, \sum_{i,k,l=1}^n (W_{ipkl} +
W_{ilkp}) \, (W_{iqkl} + W_{ilkq}) \, r^{n+5} \\
&+ \frac{1}{2n(n+2)(n+4)} \, |S^{n-1}| \, \sum_{i,j,k,l=1}^n (W_{ijkl} +
W_{ilkj})^2 \, \delta_{pq} \, r^{n+5}.
\end{align*}
\end{proposition}

\textbf{Proof.}
See \cite{Brendle1}, Proposition 16. \\

\begin{proposition}
\label{integral.identity.2}
We have
\begin{align*}
&\int_{\partial B_r(0)} \sum_{i,k,l=1}^n (\partial_l
\overline{H}_{ik}(x))^2 \, x_p \, x_q \\
&= \frac{2}{n(n+2)(n+4)} \, |S^{n-1}| \, \sum_{i,k,l=1}^n (W_{ipkl} +
W_{ilkp}) \, (W_{iqkl} + W_{ilkq}) \\
&\hspace{10mm} \cdot r^{n+3} \, \Big [ (n+4) \, f(r^2)^2 + 8r^2 \,
f(r^2) \, f'(r^2) + 4r^4 \, f'(r^2)^2 \Big ] \\
&+ \frac{1}{n(n+2)(n+4)} \, |S^{n-1}| \, \sum_{i,j,k,l=1}^n (W_{ijkl} +
W_{ilkj})^2 \, \delta_{pq} \\
&\hspace{10mm} \cdot r^{n+3} \, \Big [ (n+4) \, f(r^2)^2 + 4r^2 \,
f(r^2) \, f'(r^2) + 2r^4 \, f'(r^2)^2 \Big ].
\end{align*}
\end{proposition}

\textbf{Proof.}
Using the identity
\[\partial_l \overline{H}_{ik}(x) = f(|x|^2) \, \partial_l H_{ik}(x) + 2
\, f'(|x|^2) \, H_{ik}(x) \, x_l\]
and Euler's theorem, we obtain
\begin{align*}
&\sum_{i,k,l=1}^n (\partial_l \overline{H}_{ik}(x))^2 \\
&= f(|x|^2)^2 \, \sum_{i,k,l=1}^n (\partial_l H_{ik}(x))^2 \\
&+ 4 \, f(|x|^2) \, f'(|x|^2) \, \sum_{i,k,l=1}^n H_{ik}(x) \, x_l \, \partial_l H_{ik}(x) \\ 
&+ 4 \, |x|^2 \, f'(|x|^2)^2 \, \sum_{i,k=1}^n H_{ik}(x)^2 \\
&= f(|x|^2)^2 \, \sum_{i,k,l=1}^n (\partial_l H_{ik}(x))^2 \\
&+ \big [ 8 \, f(|x|^2) \, f'(|x|^2) + 4 \, |x|^2 \, f'(|x|^2)^2 \big ]
\, \sum_{i,k=1}^n H_{ik}(x)^2.
\end{align*}
Hence, the assertion follows from the previous proposition. \\

\begin{corollary}
\label{integral.identity.3}
We have
\begin{align*}
&\int_{\partial B_r(0)} \sum_{i,k,l=1}^n (\partial_l
\overline{H}_{ik}(x))^2 \\ &= \frac{1}{n(n+2)} \, |S^{n-1}| \,
\sum_{i,j,k,l=1}^n (W_{ijkl} + W_{ilkj})^2 \\
&\hspace{10mm} \cdot r^{n+1} \, \Big [ (n+2) \, f(r^2)^2 + 4 \, r^2 \,
f(r^2) \, f'(r^2)
+ 2 \, r^4 \, f'(r^2)^2 \Big ].
\end{align*}
\end{corollary}

\vspace{2mm}

\begin{proposition}
\label{formula.for.F}
We have
\begin{align*}
F(0,\varepsilon) &= -\frac{n-2}{16n(n-1)(n+2)} \, |S^{n-1}| \,
\sum_{i,j,k,l=1}^n (W_{ijkl} + W_{ilkj})^2 \\
&\hspace{5mm} \cdot \int_0^\infty \varepsilon^{n-2} \, (\varepsilon^2 +
r^2)^{2-n} \, r^{n+1} \\
&\hspace{15mm} \cdot \Big [ (n+2) \, f(r^2)^2 + 4 \, r^2 \, f(r^2) \,
f'(r^2)
+ 2 \, r^4 \, f'(r^2)^2 \Big ] \, dr.
\end{align*}
\end{proposition}

\textbf{Proof.}
Note that $z_{(0,\varepsilon)}(x) = 0$ for all $x \in \mathbb{R}^n$.
This implies
\[F(0,\varepsilon) = -\int_{\mathbb{R}^n} \frac{n-2}{16(n-1)} \,
\varepsilon^{n-2} \, (\varepsilon^2 + |x|^2)^{2-n} \, \sum_{i,k,l=1}^n
(\partial_l \overline{H}_{ik}(x))^2.\]
Using Corollary \ref{integral.identity.3}, we obtain
\begin{align*}
&\int_{\mathbb{R}^n} \varepsilon^{n-2} \, (\varepsilon^2 + |x|^2)^{2-n}
\, \sum_{i,k,l=1}^n (\partial_l \overline{H}_{ik}(x))^2 \\
&= \frac{1}{n(n+2)} \, |S^{n-1}| \, \sum_{i,j,k,l=1}^n (W_{ijkl} +
W_{ilkj})^2 \\
&\hspace{5mm} \cdot \int_0^\infty \varepsilon^{n-2} \, (\varepsilon^2 +
r^2)^{2-n} \, r^{n+1} \\
&\hspace{15mm} \cdot \Big [ (n+2) \, f(r^2)^2 + 4 \, r^2 \, f(r^2) \,
f'(r^2)
+ 2 \, r^4 \, f'(r^2)^2 \Big ].
\end{align*}
This proves the assertion. \\

\begin{proposition}
\label{i}
The function $F(0,\varepsilon)$ can be written in the form
\begin{align*}
F(0,\varepsilon) &= -\frac{n-2}{16n(n-1)(n+2)} \, |S^{n-1}| \,
\sum_{i,j,k,l=1}^n (W_{ijkl} + W_{ilkj})^2 \\
&\hspace{10mm} \cdot I(\varepsilon^2) \, \int_0^\infty (1 + r^2)^{2-n}
\, r^{n+7} \, dr,
\end{align*}
where
\begin{align*}
I(s) &= \frac{n-12}{n+6} \, \frac{n-10}{n+4} \, (n-8) \, \tau^2 \, s^2 +
10 \, \frac{n-12}{n+6} \, (n-10) \, \tau \, s^3 \\
&+ \Big ( 25 \, \frac{n-12}{n+6} \, (n+8) - 2(n-12) \, \tau \Big ) \,
s^4 + \Big ( \frac{n+8}{10} \, \tau - 10(n+12) \Big ) \, s^5 \\
&+ \frac{n+8}{n-14} \, \frac{3n+52}{2} \, s^6 - \frac{n+8}{n-14} \,
\frac{n+10}{n-16} \, \frac{n+24}{10} \, s^7 \\
&+ \frac{n+8}{n-14} \, \frac{n+10}{n-16} \, \frac{n+12}{n-18} \,
\frac{n+32}{400} \, s^8.
\end{align*}
\end{proposition}

\textbf{Proof.} It is straightforward to check that
\begin{align*}
&(n+2) \, f(s)^2 + 4s \, f(s) \, f'(s) + 2s^2 \, f'(s)^2 \\
&= (n+2)\tau^2 + 10(n+4)\tau \, s + \Big ( 25(n+8) - 2(n+6)\tau \Big )
\, s^2 \\
&+ \Big ( \frac{n+8}{10} \, \tau - 10(n+12) \Big ) \, s^3 +
\frac{3n+52}{2} \, s^4 - \frac{n+24}{10} \, s^5
+ \frac{n+32}{400} \, s^6.
\end{align*}
This implies
\begin{align*}
&\int_0^\infty \varepsilon^{n-2} \, (\varepsilon^2 + r^2)^{2-n} \,
r^{n+1} \\
&\hspace{10mm} \cdot \Big [ (n+2) \, f(r^2)^2 + 4 \, r^2 \, f(r^2) \,
f'(r^2)
+ 2r^4 \, f'(r^2)^2 \Big ] \, dr \\
&= (n+2)\tau^2 \, \varepsilon^4 \int_0^\infty (1+r^2)^{2-n} \, r^{n+1}
\, dr \\
&+ 10(n+4)\tau \, \varepsilon^6 \int_0^\infty (1+r^2)^{2-n} \, r^{n+3}
\, dr \\
&+ \Big ( 25(n+8) - 2(n+6)\tau \Big ) \, \varepsilon^8 \int_0^\infty
(1+r^2)^{2-n} \, r^{n+5} \, dr \\
&+ \Big ( \frac{n+8}{10} \, \tau - 10(n+12) \Big ) \, \varepsilon^{10}
\int_0^\infty (1+r^2)^{2-n} \, r^{n+7} \, dr \\
&+ \frac{3n+52}{2} \, \varepsilon^{12} \int_0^\infty (1+r^2)^{2-n} \,
r^{n+9} \, dr \\
&- \frac{n+24}{10} \, \varepsilon^{14} \int_0^\infty (1+r^2)^{2-n} \,
r^{n+11} \, dr \\
&+ \frac{n+32}{400} \, \varepsilon^{16} \int_0^\infty (1+r^2)^{2-n} \,
r^{n+13} \, dr.
\end{align*}
Using the identity
\[\int_0^\infty (1+r^2)^{2-n} \, r^{\beta+2} \, dr
= \frac{\beta+1}{2n-\beta-7} \int_0^\infty (1+r^2)^{2-n}
\, r^\beta \, dr,\]
we obtain
\begin{align*}
&\int_0^\infty \varepsilon^{n-2} \, (\varepsilon^2 + r^2)^{2-n} \,
r^{n+1} \\
&\hspace{10mm} \cdot \Big [ (n+2) \, f(r^2)^2 + 4 \, r^2 \, f(r^2) \,
f'(r^2)
+ 2r^4 \, f'(r^2)^2 \Big ] \, dr \\
&= I(\varepsilon^2) \, \int_0^\infty (1+r^2)^{2-n} \, r^{n+7} \, dr.
\end{align*}
This completes the proof. \\

In the next step, we compute the Hessian of $F$ at $(0,\varepsilon)$. \\

\begin{proposition}
\label{Hessian.of.F.1}
The second order partial derivatives of the function
$F(\xi,\varepsilon)$ are given by
\begin{align*}
\frac{\partial^2}{\partial \xi_p \, \partial \xi_q} F(0,\varepsilon)
&= \int_{\mathbb{R}^n} (n-2)^2 \, \varepsilon^{n-2} \, (\varepsilon^2 +
|x|^2)^{-n} \, \sum_{l=1}^n \overline{H}_{pl}(x) \, \overline{H}_{ql}(x) \\
&- \int_{\mathbb{R}^n} \frac{(n-2)^2}{4} \, \varepsilon^{n-2} \,
(\varepsilon^2 + |x|^2)^{-n} \, \sum_{i,k,l=1}^n (\partial_l
\overline{H}_{ik}(x))^2 \, x_p \, x_q \\
&+ \int_{\mathbb{R}^n} \frac{(n-2)^2}{8(n-1)} \, \varepsilon^{n-2} \,
(\varepsilon^2 + |x|^2)^{1-n} \, \sum_{i,k,l=1}^n (\partial_l
\overline{H}_{ik}(x))^2 \, \delta_{pq}.
\end{align*}
\end{proposition}

\textbf{Proof.}
See \cite{Brendle1}, Proposition 21. \\

\begin{proposition}
\label{Hessian.of.F.2}
The second order partial derivatives of the function
$F(\xi,\varepsilon)$ are given by
\begin{align*}
&\frac{\partial^2}{\partial \xi_p \, \partial \xi_q} F(0,\varepsilon) \\
&= -\frac{2(n-2)^2}{n(n+2)(n+4)} \, |S^{n-1}| \, \sum_{i,k,l=1}^n
(W_{ipkl} + W_{ilkp}) \, (W_{iqkl} + W_{ilkq}) \\
&\hspace{5mm} \cdot \int_0^\infty \varepsilon^{n-2} \, (\varepsilon^2 +
r^2)^{-n} \, r^{n+5} \, \Big [ 2 \, f(r^2) \, f'(r^2) + r^2 \, f'(r^2)^2
\Big ] \, dr \\
&- \frac{(n-2)^2}{2n(n+2)(n+4)} \, |S^{n-1}| \, \sum_{i,j,k,l=1}^n
(W_{ijkl} + W_{ilkj})^2 \, \delta_{pq} \\
&\hspace{5mm} \cdot \int_0^\infty \varepsilon^{n-2} \, (\varepsilon^2 +
r^2)^{-n} \, r^{n+5} \, \Big [
2 \, f(r^2) \, f'(r^2) + r^2 \, f'(r^2)^2 \Big ] \, dr \\
&+ \frac{(n-2)^2}{4n(n-1)(n+2)} \, |S^{n-1}| \, \sum_{i,j,k,l=1}^n
(W_{ijkl} + W_{ilkj})^2 \, \delta_{pq} \\
&\hspace{5mm} \cdot \int_0^\infty \varepsilon^{n-2} \, (\varepsilon^2 +
r^2)^{1-n} \, r^{n+5} \, f'(r^2)^2 \, dr.
\end{align*}
\end{proposition}

\textbf{Proof.}
Using the identity
\begin{align*}
&\int_{\partial B_r(0)} \sum_{l=1}^n \overline{H}_{pl}(x) \,
\overline{H}_{ql}(x) \\
&= \int_{\partial B_r(0)} \sum_{i,j,k,l,m=1}^n W_{ipkl} \, W_{jqml} \,
x_i \, x_j \, x_k \, x_m \, f(|x|^2)^2 \\
&= \frac{1}{n(n+2)} \, |S^{n-1}| \\
&\hspace{10mm} \cdot \sum_{i,j,k,l,m=1}^n W_{ipkl} \, W_{jqml} \,
(\delta_{ij} \, \delta_{km} + \delta_{ik} \, \delta_{jm} + \delta_{im}
\, \delta_{jk}) \, r^{n+3} \, f(r^2)^2 \\
&= \frac{1}{2n(n+2)} \, |S^{n-1}| \, \sum_{i,k,l=1}^n (W_{ipkl} +
W_{ilkp}) \, (W_{iqkl} + W_{ilkq}) \, r^{n+3} \, f(r^2)^2,
\end{align*}
we obtain
\begin{align*}
&\int_{\mathbb{R}^n} \varepsilon^{n-2}
\, (\varepsilon^2 + |x|^2)^{-n} \, \sum_{i,k,l=1}^n \overline{H}_{pl}(x)
\, \overline{H}_{ql}(x) \\
&= \frac{1}{2n(n+2)} \, |S^{n-1}| \,
\sum_{i,k,l=1}^n (W_{ipkl} + W_{ilkp}) \, (W_{iqkl} + W_{ilkq}) \\
&\hspace{5mm} \cdot \int_0^\infty
\varepsilon^{n-2} \, (\varepsilon^2 + r^2)^{-n} \, r^{n+3} \, f(r^2)^2
\, dr.
\end{align*}
Similarly, it follows from Proposition \ref{integral.identity.2} that
\begin{align*}
&\int_{\mathbb{R}^n} \varepsilon^{n-2} \, (\varepsilon^2 + |x|^2)^{-n}
\, \sum_{i,k,l=1}^n (\partial_l \overline{H}_{ik}(x))^2 \, x_p \, x_q \\
&= \frac{2}{n(n+2)(n+4)} \, |S^{n-1}| \, \sum_{i,k,l=1}^n (W_{ipkl} +
W_{ilkp}) \, (W_{iqkl} + W_{ilkq}) \\
&\hspace{5mm} \cdot \int_0^\infty \varepsilon^{n-2} \, (\varepsilon^2 +
r^2)^{-n} \, r^{n+3} \\ 
&\hspace{15mm} \cdot \Big [ (n+4) \, f(r^2)^2 + 8r^2 \, f(r^2) \,
f'(r^2) + 4r^4 \, f'(r^2)^2 \Big ] \, dr \\
&+ \frac{1}{n(n+2)(n+4)} \, |S^{n-1}| \, \sum_{i,j,k,l=1}^n (W_{ijkl} +
W_{ilkj})^2 \, \delta_{pq} \\
&\hspace{5mm} \cdot \int_0^\infty \varepsilon^{n-2} \, (\varepsilon^2 +
r^2)^{-n} \, r^{n+3} \\ 
&\hspace{15mm} \cdot \Big [ (n+4) \, f(r^2)^2 + 4r^2 \, f(r^2) \,
f'(r^2) + 2r^4 \, f'(r^2)^2 \Big ] \, dr.
\end{align*}
Moreover, we have
\begin{align*}
&\int_{\mathbb{R}^n} \varepsilon^{n-2} \, (\varepsilon^2 + |x|^2)^{1-n}
\, \sum_{i,k,l=1}^n (\partial_l \overline{H}_{ik}(x))^2 \, \delta_{pq} \\
&= \frac{1}{n(n+2)} \, |S^{n-1}| \, \sum_{i,j,k,l=1}^n (W_{ijkl} +
W_{ilkj})^2 \, \delta_{pq} \\
&\hspace{5mm} \cdot \int_0^\infty \varepsilon^{n-2} \, (\varepsilon^2 +
r^2)^{1-n} \, r^{n+1} \\ 
&\hspace{15mm} \cdot \Big [ (n+2) \, f(r^2)^2 + 4 \, r^2 \, f(r^2)
\, f'(r^2) + 2 \, r^4 \, f'(r^2)^2 \Big ] \, dr
\end{align*}
by Corollary \ref{integral.identity.3}. A straightforward calculation yields
\begin{align*}
&(\varepsilon^2+r^2)^{1-n} \, r^{n+1} \, \big [ (n+2) \, f(r^2)^2 + 4r^2
\, f(r^2) \, f'(r^2) \big ] \\
&= 2(n-1) \, (\varepsilon^2+r^2)^{-n} \, r^{n+3} \, f(r^2)^2 +
\frac{d}{dr} \big [ (\varepsilon^2+r^2)^{1-n} \, r^{n+2} \, f(r^2)^2 \big ].
\end{align*}
This implies
\begin{align*}
&\int_{\mathbb{R}^n} \varepsilon^{n-2} \, (\varepsilon^2 + |x|^2)^{1-n}
\, \sum_{i,k,l=1}^n (\partial_l \overline{H}_{ik}(x))^2 \, \delta_{pq} \\
&= \frac{2(n-1)}{n(n+2)} \, |S^{n-1}| \, \sum_{i,j,k,l=1}^n (W_{ijkl} +
W_{ilkj})^2 \, \delta_{pq} \\ 
&\hspace{5mm} \cdot
\int_0^\infty \varepsilon^{n-2} \, (\varepsilon^2 + r^2)^{-n} \, r^{n+3}
\, f(r^2)^2 \, dr \\
&+ \frac{2}{n(n+2)} \, |S^{n-1}| \, \sum_{i,j,k,l=1}^n (W_{ijkl} +
W_{ilkj})^2 \, \delta_{pq} \\
&\hspace{5mm} \cdot \int_0^\infty \varepsilon^{n-2} \, (\varepsilon^2 +
r^2)^{1-n} \, r^{n+5} \, f'(r^2)^2 \, dr.
\end{align*}
Putting these facts together, the assertion follows. \\

\begin{proposition}
\label{j}
We have
\begin{align*}
&\int_0^\infty \varepsilon^{n-2} \, (\varepsilon^2 + r^2)^{-n} \,
r^{n+5} \, \Big [ 2 \, f(r^2) \, f'(r^2) + r^2 \, f'(r^2)^2 \Big ] \, dr \\
&= J(\varepsilon^2) \, \int_0^\infty (1+r^2)^{-n} \, r^{n+9} \, dr,
\end{align*}
where
\begin{align*}
J(s) &= 10 \, \frac{n-10}{n+8} \, \frac{n-8}{n+6} \, \tau \, s^2 +
\frac{n-10}{n+8} \, (75-4\tau) \, s^3 \\
&+ \Big ( \frac{3}{10} \, \tau - 50 \Big ) \, s^4 + \frac{23}{2} \,
\frac{n+10}{n-12} \, s^5 - \frac{11}{10} \, \frac{n+10}{n-12} \,
\frac{n+12}{n-14} \, s^6 \\
&+ \frac{3}{80} \, \frac{n+10}{n-12} \, \frac{n+12}{n-14} \,
\frac{n+14}{n-16} \, s^7.
\end{align*}
\end{proposition}

\textbf{Proof.}
Note that
\begin{align*}
&2 \, f(s) \, f'(s) + s \, f'(s)^2 \\
&= 10\tau + (75-4\tau) \, s + \Big ( \frac{3}{10} \, \tau - 50 \Big ) \,
s^2 + \frac{23}{2} \, s^3
- \frac{11}{10} \, s^4 + \frac{3}{80} \, s^5.
\end{align*}
This implies
\begin{align*}
&\int_0^\infty \varepsilon^{n-2} \, (\varepsilon^2 + r^2)^{-n} \,
r^{n+5} \, \Big [ 2 \, f(r^2) \, f'(r^2)
+ r^2 \, f'(r^2)^2 \Big ] \, dr \\
&= 10\tau \, \varepsilon^4 \, \int_0^\infty (1 + r^2)^{-n} \, r^{n+5} \,
dr \\
&+ (75-4\tau) \, \varepsilon^6 \, \int_0^\infty (1 + r^2)^{-n} \,
r^{n+7} \, dr \\
&+ \Big ( \frac{3}{10} \, \tau - 50 \Big ) \, \varepsilon^8 \,
\int_0^\infty (1 + r^2)^{-n} \, r^{n+9} \, dr \\
&+ \frac{23}{2} \, \varepsilon^{10} \, \int_0^\infty (1 + r^2)^{-n} \,
r^{n+11} \, dr \\
&- \frac{11}{10} \, \varepsilon^{12} \, \int_0^\infty (1 + r^2)^{-n} \,
r^{n+13} \, dr \\
&+ \frac{3}{80} \, \varepsilon^{14} \, \int_0^\infty (1 + r^2)^{-n} \,
r^{n+15} \, dr.
\end{align*}
Hence, the assertion follows from the identity
\[\int_0^\infty (1+r^2)^{-n} \, r^{\beta+2} \, dr
= \frac{\beta+1}{2n-\beta-3} \int_0^\infty (1+r^2)^{-n}
\, r^\beta \, dr.\]

\begin{proposition}
Assume that $25 \leq n \leq 51$. Then we can choose $\tau \in
\mathbb{R}$ such that $I'(1) = 0$, $I''(1) < 0$, and $J(1) < 0$.
\end{proposition}

\textbf{Proof.}
The condition $I'(1) = 0$ is equivalent to
\[a_n \, \tau^2 + b_n \, \tau + c_n = 0,\]
where
\begin{align*}  
a_n &= 2 \, \frac{n-12}{n+6} \, \frac{n-10}{n+4} \, (n-8) \\
b_n &= 30 \, \frac{n-12}{n+6} \, (n-10) - 8(n-12) + \frac{n+8}{2} \\
c_n &= 100 \, \frac{n-12}{n+6} \, (n+8) - 50(n+12) + 3 \,
\frac{n+8}{n-14} \, (3n+52) \\ &- 7 \, \frac{n+8}{n-14} \,
\frac{n+10}{n-16} \, \frac{n+24}{10} + \frac{n+8}{n-14} \,
\frac{n+10}{n-16} \, \frac{n+12}{n-18} \, \frac{n+32}{50}.
\end{align*}
By inspection, one verifies that $49 \, a_n - 7 \, b_n + c_n < 0$ for
$25 \leq n \leq 51$. Since $a_n$ is positive, there exists a unique real
number $\tau < -7$ such that $a_n \, \tau^2 + b_n \, \tau + c_n = 0$.
Moreover, we have
\[I''(1) - I'(1) = \alpha_n \, \tau + \beta_n\]
and
\[J(1) = \gamma_n \, \tau + \delta_n,\]
where
\begin{align*}
\alpha_n &= 30 \, \frac{n-12}{n+6} \, (n-10) - 16(n-12) +
\frac{3(n+8)}{2} \\
\beta_n &= 200 \, \frac{n-12}{n+6} \, (n+8) - 150(n+12) + 12 \,
\frac{n+8}{n-14} \, (3n+52) \\ &- 35 \, \frac{n+8}{n-14} \,
\frac{n+10}{n-16} \, \frac{n+24}{10} + 3 \, \frac{n+8}{n-14} \,
\frac{n+10}{n-16} \, \frac{n+12}{n-18} \, \frac{n+32}{25} \\[2mm]
\gamma_n &= 10 \, \frac{n-10}{n+8} \, \frac{n-8}{n+6} -
\frac{4(n-10)}{n+8} + \frac{3}{10} \\
\delta_n &= 75 \, \frac{n-10}{n+8} - 50 + \frac{23}{2} \,
\frac{n+10}{n-12} - \frac{11}{10} \, \frac{n+10}{n-12} \,
\frac{n+12}{n-14} \\ &+ \frac{3}{80} \, \frac{n+10}{n-12} \,
\frac{n+12}{n-14} \, \frac{n+14}{n-16}.
\end{align*}
By inspection, one verifies that $7\alpha_n > \beta_n > 0$ and
$7\gamma_n > \delta_n > 0$ for $25 \leq n \leq 51$. This implies
$I''(1) = \alpha_n \, \tau + \beta_n < -7\alpha_n + \beta_n < 0$ and
$J(1) = \gamma_n \, \tau + \delta_n < -7\gamma_n + \delta_n < 0$. This
completes the proof. \\

\begin{corollary}
\label{strict.local.minimum}
Assume that $\tau$ is chosen such that $I'(1) = 0$, $I''(1) < 0$, and
$J(1) < 0$. Then the function $F(\xi,\varepsilon)$ has a strict local
minimum at $(0,1)$.
\end{corollary}

\textbf{Proof.}
Since $I'(1) = 0$, we have $\frac{\partial}{\partial \varepsilon} F(0,1)
= 0$. Therefore, $(0,1)$ is a critical point of the function
$F(\xi,\varepsilon)$. Since $J(1) < 0$, we have
\[\int_0^\infty (1 + r^2)^{-n} \, r^{n+5} \, \Big [ 2 \, f(r^2) \,
f'(r^2) + r^2 \, f'(r^2)^2 \Big ] \, dr < 0\]
by Proposition \ref{j}. Hence, it follows from Proposition
\ref{Hessian.of.F.2} that the matrix $\frac{\partial^2}{\partial \xi_p
\, \partial \xi_q} F(0,1)$ is positive definite. Using 
Proposition \ref{i} and the inequality $I''(0) < 0$, we obtain 
$\frac{\partial^2}{\partial \varepsilon^2} F(0,1) > 0$. Consequently, 
the function $F(\xi,\varepsilon)$ has a strict local minimum at $(0,1)$. \\

\section{Proof of the main theorem}

\begin{proposition}
\label{perturbation.argument}
Assume that $25 \leq n \leq 51$. Moreover, let $g$ be a smooth metric on
$\mathbb{R}^n$ of the form
$g(x) = \exp(h(x))$, where $h(x)$ is a trace-free symmetric two-tensor
on $\mathbb{R}^n$ such
that $|h(x)| + |\partial h(x)| + |\partial^2 h(x)| \leq \alpha \leq
\alpha_1$ for all $x \in \mathbb{R}^n$,
$h(x) = 0$ for $|x| \geq 1$, and
\[h_{ik}(x) = \mu \, \lambda^6 \, f(\lambda^{-2} \, |x|^2) \,
H_{ik}(x)\] for $|x| \leq \rho$. If $\alpha$ and $\rho^{2-n} \, \mu^{-2}
\, \lambda^{n-18}$ are
sufficiently small, then there exists a positive function $v$ such that
\[\Delta_g v - \frac{n-2}{4(n-1)} \, R_g \, v + n(n-2) \,
v^{\frac{n+2}{n-2}} = 0,\]
\[\int_{\mathbb{R}^n} v^{\frac{2n}{n-2}} < \Big ( \frac{Y(S^n)}{4n(n-1)}
\Big )^{\frac{n}{2}},\] and  
$\sup_{|x| \leq \lambda} v(x) \geq c \, \lambda^{\frac{2-n}{2}}$. Here,
$c$ is a positive constant that depends only on $n$.
\end{proposition}

\textbf{Proof.}
By Corollary \ref{strict.local.minimum}, the function
$F(\xi,\varepsilon)$ has a strict local minimum at $(0,1)$. It follows
from Proposition
\ref{formula.for.F} that $F(0,1) < 0$. Hence, we can find an open set
$\Omega' \subset \Omega$ such that $(0,1) \in \Omega'$ and
\[F(0,1) < \inf_{(\xi,\varepsilon) \in \partial \Omega'}
F(\xi,\varepsilon) < 0.\]
Using Corollary \ref{key.estimate}, we obtain
\begin{align*}
&|\mathcal{F}_g(\lambda\xi,\lambda\varepsilon) - \lambda^{16} \, \mu^2
\, F(\xi,\varepsilon)| \\
&\leq C \, \lambda^{\frac{16n}{n-2}} \, \mu^{\frac{2n}{n-2}} + C \,
\lambda^8 \, \mu \, \Big ( \frac{\lambda}{\rho} \Big )^{\frac{n-2}{2}} +
C \, \Big ( \frac{\lambda}{\rho} \Big )^{n-2}
\end{align*}
for all $(\xi,\varepsilon) \in \Omega$. This implies
\begin{align*}
&|\lambda^{-16} \, \mu^{-2} \,
\mathcal{F}_g(\lambda\xi,\lambda\varepsilon) - F(\xi,\varepsilon)| \\
&\leq C \, \lambda^{\frac{32}{n-2}} \, \mu^{\frac{4}{n-2}} + C \,
\rho^{\frac{2-n}{2}} \, \mu^{-1} \, \lambda^{\frac{n-18}{2}} + C \,
\rho^{2-n} \, \mu^{-2} \, \lambda^{n-18}
\end{align*}
for all $(\xi,\varepsilon) \in \Omega$.
Hence, if $\rho^{2-n} \, \mu^{-2} \, \lambda^{n-18}$ is sufficiently
small, then we have
\[\mathcal{F}_g(0,\lambda) < \inf_{(\xi,\varepsilon) \in \partial
\Omega'} \mathcal{F}_g(\lambda\xi,\lambda\varepsilon) < 0.\]
Consequently, there exists a point $(\bar{\xi},\bar{\varepsilon}) \in
\Omega'$ such that
\[\mathcal{F}_g(\lambda\bar{\xi},\lambda\bar{\varepsilon}) =
\inf_{(\xi,\varepsilon) \in \Omega'}
\mathcal{F}_g(\lambda\xi,\lambda\varepsilon) < 0.\]
By Proposition \ref{reduction.to.a.finite.dimensional.problem}, the
function $v = v_{(\lambda\bar{\xi},\lambda\bar{\varepsilon})}$ is a
non-negative weak solution of the partial differential equation
\[\Delta_g v - \frac{n-2}{4(n-1)} \, R_g \, v + n(n-2) \,
v^{\frac{n+2}{n-2}} = 0.\]
Using a result of N.~Trudinger, we conclude that $v$ is smooth (see
\cite{Trudinger}, Theorem 3 on p. 271). Moreover, we have
\begin{align*}
2(n-2) \int_{\mathbb{R}^n} v^{\frac{2n}{n-2}}
&= 2(n-2) \, \Big ( \frac{Y(S^n)}{4n(n-1)} \Big )^{\frac{n}{2}} +
\mathcal{F}_g(\lambda\bar{\xi},\lambda\bar{\varepsilon}) \\
&< 2(n-2) \, \Big ( \frac{Y(S^n)}{4n(n-1)} \Big )^{\frac{n}{2}}.
\end{align*}
Finally, it follows from Proposition \ref{fixed.point.argument} that
$\|v -
u_{(\lambda\bar{\xi},\lambda\bar{\varepsilon})}\|_{L^{\frac{2n}{n-2}}(\mathbb{R}^n)}
\leq C \, \alpha$.
This implies
\[|B_\lambda(0)|^{\frac{n-2}{2n}} \, \sup_{|x| \leq \lambda} v(x) \geq
\|v\|_{L^{\frac{2n}{n-2}}(B_\lambda(0))} \geq
\|u_{(\lambda\bar{\xi},\lambda\bar{\varepsilon})}\|_{L^{\frac{2n}{n-2}}(B_\lambda(0))}
- C \, \alpha.\]
Hence, if $\alpha$ is sufficiently small, then we obtain
$\lambda^{\frac{n-2}{2}} \, \sup_{|x| \leq \lambda} v(x) \geq c$. \\

\begin{proposition}
Let $25 \leq n \leq 51$. Then there exists a smooth metric $g$ on
$\mathbb{R}^n$ with the following properties:
\begin{itemize}
\item[(i)] $g_{ik}(x) = \delta_{ik}$ for $|x| \geq \frac{1}{2}$
\item[(ii)] $g$ is not conformally flat
\item[(iii)] There exists a sequence of non-negative smooth functions
$v_\nu$ ($\nu \in \mathbb{N}$) such that  
\[\Delta_g v_\nu - \frac{n-2}{4(n-1)} \, R_g \, v_\nu + n(n-2) \,
v_\nu^{\frac{n+2}{n-2}} = 0\]
for all $\nu \in \mathbb{N}$,
\[\int_{\mathbb{R}^n} v_\nu^{\frac{2n}{n-2}} < \Big (
\frac{Y(S^n)}{4n(n-1)} \Big )^{\frac{n}{2}}\]
for all $\nu \in \mathbb{N}$, and $\sup_{|x| \leq 1} v_\nu(x) \to
\infty$ as $\nu \to \infty$.
\end{itemize}
\end{proposition}

\textbf{Proof.}
Choose a smooth cutoff function $\eta: \mathbb{R} \to \mathbb{R}$ such
that $\eta(t) = 1$ for $t \leq 1$ and $\eta(t) = 0$ for $t \geq 2$. We
 define a trace-free symmetric two-tensor on $\mathbb{R}^n$ by
\[h_{ik}(x) = \sum_{N=N_0}^\infty \eta(4N^2 \, |x - y_N|) \, 2^{-4N} \,
f(2^{N} \, |x - y_N|^2) \, H_{ik}(x - y_N),\]
where $y_N = (\frac{1}{N},0,\hdots,0) \in \mathbb{R}^n$. It is
straightforward to verify that $h(x)$ is $C^\infty$ smooth.
Moreover, if $N_0$ is sufficiently large, then we have $h(x) = 0$ for
$|x| \geq \frac{1}{2}$ and $|h(x)| + |\partial h(x)| + |\partial^2 h(x)|
\leq \alpha$ for all $x \in \mathbb{R}^n$.
(Here, $\alpha$ is the constant that appears in Proposition
\ref{perturbation.argument}.) We now define a Riemannian metric $g$ by
$g(x) = \exp(h(x))$. The assertion is then a consequence of Proposition
\ref{perturbation.argument}. \\

\end{document}